# On a Theorem of Lenstra and Schoof

B.V. Petrenko [1]

*Department of Mathematics, University of Illinois, 1409 West Green Street, Urbana, Illinois 61801, U.S.A.*

**Abstract**

We give a detailed proof of Theorem 1.15 from a well-known paper "Primitive normal bases for finite fields" by H.W. Lenstra Jr. and R.J. Schoof. We are not aware of any other proofs. Let $L/K$ be a finite-dimensional Galois field extension and $B$ the set of all normal bases of this extension. Theorem 1.15 describes the group of all $\gamma$ in the multiplicative group of $L$ such that $\gamma B = B$.

## 1 Introduction

The celebrated paper of Lenstra and Schoof [6] solved a difficult problem. It proved that every extension of finite fields has a primitive normal basis. This means that for any extension of finite fields $\mathbb{F}_{q^m}/\mathbb{F}_q$, there exists $\alpha \in \mathbb{F}_{q^m}$ such that $\alpha$ generates the multiplicative group $\mathbb{F}_{q^m}^*$ of $\mathbb{F}_{q^m}$, and $\{\alpha^{q^j} \mid j \in \{1, ..., m\}\}$ is an $\mathbb{F}_q$-basis of $\mathbb{F}_{q^m}$. Several particular cases have been done earlier by Carlitz in [1] and [2], and Davenport in [5]. Statement 1.15 is a step in the proof of the main result of [6]. We recall this statement below for the convenience of the reader. We call it Theorem 1 for reference purposes.

*Email address:* petrenko@uiuc.edu (B.V. Petrenko).
[1] I thank Hendrik Lenstra and Everett Dade for their ideas and help.

**Theorem 1** *Let $K \subseteq L$ be a finite extension of fields with Galois group $G$. Let*
$$B = \{\alpha \in L \mid (\tau(\alpha))_{\tau \in G} \text{ is a basis of } L \text{ over } K\},$$
*and denote by $w$ the number $|G|$th roots of unity in $K^*$. Then for $\gamma \in L^*$ the following assertions are equivalent:*

*(i)* $\gamma B \subseteq B$.

*(ii)* $\gamma B = B$.

*(iii)* $\frac{\tau(\gamma)}{\gamma} \in K^*$ *for all* $\tau \in G$.

*(iv)* $\gamma^w \in K^*$.

*The set of all $\gamma \in L^*$ satisfying these conditions is a subgroup of $L^*$ containing $K^*$ and $C/K^*$ is isomorphic to the group of all homomorphisms $G \to K^*$.*

This theorem is stated in [6] without a proof and used there in the case when $L$ is a finite field. Because we are not aware of any proof of this result, we prove it below using the Wedderburn's Theorems (see, for example, [3]). This is the main result of the present paper. It is established in Section 4. Our proof is based on a private communication of Lenstra [6].

## 2 Notation

If $W$ is a set, we denote its cardinality by $|W|$.

$F$ is a field

$F^*$ is the multiplicative group of $F$

$A$ is a finite-dimensional associative algebra over $F$ with a two-sided identity $1_A$. We assume that $0_A \neq 1_A$.

$J(A)$ is the Jacobson radical of $A$.

$U(A)$ is the unit group of $A$, i.e. the set of $u \in A$ having two-sided multiplicative inverses $u^{-1}$.

$\operatorname{nil}(A)$ the set of all nilpotent elements of $A$.

$H$ is a hyperplane in $A$, i.e. an $F$-subspace of $A$ of codimension one.

$\mathbb{F}_2$ is a field with two elements.

$M_n(R)$ is the ring of $n \times n$ matrices with entries in a ring $R$.

$E_{i,j}$ is the matrix in $M_n(R)$ all of whose entries are $0_R$ except the $(i,j)$-entry which is $1_R$. It will always be clear what $R$ and $n$ are in any particular situation.

$\omega FG$ is the augmentation ideal of the group algebra $FG$, where $F$ is a field and $G$ is a group.

Tr is the map $\text{Tr}: L \to K, t \mapsto \sum_{\tau \in G} \tau(t)$.

$N$ is the kernel of the map Tr, defined above.

$\mathbb{N}$ is the set of positive integers.

## 3 Preliminary Results

We begin by stating three main theorems of this section. Their proofs are based on a private communication of Dade [4].

**Theorem 2** *If $H$ is a hyperplane in an $F$-algebra $A$, then $H \cap U(A) \neq \emptyset$ except possibly when either*

*(i) $H$ is a two-sided ideal of $A$, or*

*(ii) $|F| = 2$ and $A$ has a factor algebra isomorphic to $\mathbb{F}_2 \oplus \mathbb{F}_2 \oplus \mathbb{F}_2$.*

**Theorem 3** *If $H$ is a hyperplane in an $F$-algebra $A$, then $H \setminus U(A) \neq \emptyset$ except possibly when*

*(iii) $|F| = 2$ and $A$ has a factor algebra isomorphic to $\mathbb{F}_2 \oplus \mathbb{F}_2$.*

**Theorem 4** *If $H$ is a hyperplane in the group algebra $FG$ of a finite group $G$ over a field $F$, then $H \cap U(A) \neq \emptyset$ and $H \setminus U(A) \neq \emptyset$ except when $H$ is a two-sided ideal of $A$.*

We see that Theorems 2 and 3 imply Theorem 4 except possibly when $|F| = 2$. Lemma 5 below shows that the group algebra $\mathbb{F}_2 G$ of any group $G$ has only one two-sided ideal of codimension one, namely $\omega \mathbb{F}_2 G$. Therefore, Theorems 2 and 3 imply Theorem 4 when $|F| = 2$ also.

**Lemma 5** *The group algebra $\mathbb{F}_2 G$ of any group $G$ has only one two-sided ideal of codimension one, namely $\omega \mathbb{F}_2 G$.*

**PROOF.** Let $I$ be a two-sided ideal of $\mathbb{F}_2 G$ with $\dim_{\mathbb{F}_2} \mathbb{F}_2 G / I = 1$. Then

$$\mathbb{F}_2 G \xrightarrow{\pi, \text{onto}} \mathbb{F}_2 G / I \xrightarrow{\varepsilon, \cong} \mathbb{F}_2$$

where $\pi(t) = t + I, t \in \mathbb{F}_2 G$, and $\varepsilon$ is the unique field isomorphism. We see that $G \subseteq U(\mathbb{F}_2 G)$ and $\varepsilon \circ \pi(U(\mathbb{F}_2 G)) = \{1\}$. Therefore, $\pi(g) = \pi(1)$ for any $g \in G$, i.e. $g - 1 \in I$. Since $\omega \mathbb{F}_2 G = \sum_{g \in G} \mathbb{F}_2 (g - 1)$, we have that $\omega \mathbb{F}_2 G \subseteq I$. We know that $\omega \mathbb{F}_2 G$ is a maximal two-sided ideal of $\mathbb{F}_2 G$. Therefore, $I = \omega \mathbb{F}_2 G$.

A key step in proving Theorems 2 and 3 is the following

**Lemma 6** *Let $H$ be a hyperplane in an $F$-algebra $A$. We do not assume that $\dim_F A$ is finite in this lemma. Suppose that either $U(A) \subseteq H$ or $U(A) \cap H = \emptyset$, then $\text{nil}(A) \subseteq H$.*

**PROOF.**

(1) Suppose that $U(A) \subseteq H$. Then $1_A \in H$. If $t \in \text{nil}(A)$, then $1_A + t \in U(A) \subseteq H$. Hence $t = (1_A + t) - 1_A \in H$.
(2) Suppose that $U(A) \cap H = \emptyset$. Then $1_A \notin H$. Suppose that there exists $t \in \text{nil}(A) \setminus H$. Since $\dim_F A / H = 1$, we see that $1_A + H$ and $t + H$ are $F$-linearly dependent. Therefore, $1_A + \alpha t \in H$ for some $\alpha \in F$. At the same time $\alpha t \in \text{nil}(A)$, so that $1_A + \alpha t \in U(A)$. We conclude that $U(A) \cap H \neq \emptyset$, a contradiction.

Next we prove Theorem 2 for simple algebras.

**Lemma 7** *Let $H$ be a hyperplane in a simple $F$-algebra $A$ such that $H \cap U(A) = \emptyset$. Then $A \cong F$ and $H = \{0\}$.*

**PROOF.** By Wedderburn's Theorems, $A \cong M_n(D)$ for some division $F$-algebra $D$.

If $n = 1$, then $A \cong D$, so that $U(A) = A \setminus \{0\}$. Since $H \cap U(A) = \emptyset$, we see that $H = \{0\}$. Then $\dim_F A / H = 1$ implies that $A \cong F$.

Next we show that $n \geq 2$ cannot occur. Suppose not. Let $\varphi : A \to M_n(D)$ be an isomorphism of $F$-algebras. By Lemma 6, $\text{nil}(M_n(D)) \subseteq \varphi(H)$. Hence $E_{i,j} \in \varphi(H)$ for all $i \neq j$. Therefore,

$$M = E_{n,1} + E_{1,2} + ... + E_{n-1,n} = \begin{pmatrix} 0 & 1 & 0 & \ldots & 0 \\ 0 & 0 & 1 & \ldots & 0 \\ \vdots & \vdots & \ldots & \ddots & \vdots \\ 0 & \ldots & \ldots & 0 & 1 \\ 1 & 0 & \ldots & \ldots & 0 \end{pmatrix} \in \varphi(H).$$

However, $M$ is a permutation matrix and therefore is invertible in $M_n(F) \subseteq M_n(D)$. Since $\varphi^{-1}(M) \subseteq H$, we conclude that $H \cap U(A) \neq \emptyset$, contradicting initial assumption.

Theorem 3 for simple algebras is

**Lemma 8** *No hyperplane in any simple $F$-algebra $A$ can contain $U(A)$.*

**PROOF.** Suppose not. Let $A$ be a simple $F$-algebra and $H$ a hyperplane in $A$ such that $U(A) \subseteq H$. Let $\varphi : A \to M_n(D)$ be an isomorphism of $F$-algebras, for some $n \in \mathbb{N}$ and a division $F$-algebra $D$.

If $n = 1$, then $A \cong D$ and $A \setminus \{0\} = U(A) \subseteq H$. This is impossible.

Let us consider the case $n \geq 2$. By Lemma 6, $\mathrm{nil}(M_n(D)) \subseteq \varphi(H)$. Therefore, $E_{i,j} \in \varphi(H)$ for all $i \neq j$. We also note that

$$E_{i,i} + E_{i,j} - E_{j,i} - E_{j,j} \in \varphi(H)$$

because the square of this matrix is zero. Thus

$$E_{i,i} - E_{j,j} = (E_{i,i} + E_{i,j} - E_{j,i} - E_{j,j}) - E_{i,j} + E_{j,i} \in \varphi(H).$$

Define $tr : M_n(D) \to D, (d_{i,j}) \mapsto \sum_{l=1}^{n} d_{l,l}$. Then the above discussion shows that $\mathrm{Ker}(tr) \subseteq \varphi(H)$. Hence, $\varphi(H) = tr^{-1}(H_1)$ for some $F$-hyperplane $H_1$ of $D$. Let $d \in D \setminus H_1$. Then the matrix

$$Q = (dE_{1,1} + E_{n,1}) + E_{1,2} + E_{2,3} + ... E_{n-1,n} = \begin{pmatrix} d & 1 & 0 & \ldots & 0 \\ 0 & 0 & 1 & \ldots & 0 \\ \vdots & \vdots & \ddots & \ddots & \vdots \\ 0 & \ldots & \ldots & 0 & 1 \\ 1 & 0 & \ldots & \ldots & 0 \end{pmatrix}$$

has the two-sided inverse

$$Q^{-1} = (E_{1,2} - dE_{2,n}) + E_{2,1} + E_{3,2} + ... + E_{n,n-1} = \begin{pmatrix} 0 & 0 & 0 & \ldots & \ldots & 1 \\ 1 & 0 & 0 & \ldots & \ldots & -d \\ 0 & 1 & 0 & \ldots & \ldots & 0 \\ \vdots & \vdots & \ddots & \ddots & \ldots & \vdots \\ 0 & \ldots & \ldots & \ddots & 0 & 0 \\ 0 & 0 & \ldots & \ldots & 1 & 0 \end{pmatrix}$$

and $tr(Q) = d \notin H_1$, so that $Q \in \varphi(U(A)\setminus H)$. Hence $\varphi^{-1}(Q) \in U(A)\setminus H$, a contradiction.

Now we are ready to prove Theorems 2 and 3. They hold when $\dim_F A = 1$. Below we prove these theorems simultaneously by induction on $\dim_F A$.

**PROOF OF THEOREMS 2 AND 3.** We can assume that $\dim_F A \geq 2$ and that Theorems 2 and 3 hold for all strictly smaller values of $\dim_F A$.

Assume that Theorem 2 is false for some hyperplane $H$ in $A$. Then $H \cap U(A) = \emptyset$, $H$ is not a two-sided ideal of $A$, and $A$ does not satisfy (ii).

We claim that $J(A) = \{0\}$. Suppose not. By Lemma 6, $\mathrm{nil}(A) \subseteq H$. So $J(A) \subseteq H$. Define

$$\pi : A \to A' = A/J(A), a \mapsto a + J(A).$$

Put $H' = \pi(H)$. Since $J(A) \subseteq H$, we see that that $H = \pi^{-1}(H')$ and $H'$ is a hyperplane in $A'$. Theorem 2 holds for $A'$ and $H'$. Let $u' \in H' \cap U(A')$. Then $u' = u + J(A)$ for some $u \in H$. At the same time, $u \in U(A)$ because $J(A)$ is nilpotent. Therefore $H \cap U(A) \neq \emptyset$, a contradiction.

By Wedderburn's Theorems, $A = A_1 \dotplus ... \dotplus A_n$, where $A_i$ are two-sided ideals each of which is a simple $F$-subalgebra of $A$, $n \in \mathbb{N}$.

If $n = 1$, then $A = A_1$ is simple. Since $H \cap U(A) = \emptyset$, Lemma 7 implies that $H = \{0\}$ and therefore $H$ is a two-sided ideal of $A$. This contradicts our assumptions. Hence $n \geq 2$.

We claim that $A_i \not\subseteq H$ for all $i \in \{1, ..., n\}$. Suppose not. Then $A_1 \subseteq H$ without loss of generality. Define $A' = \sum_{j=2}^n A_j$, $H' = H \cap A'$ and

$$\pi: A \to A', a_1 + \ldots a_n \mapsto a_2 + \ldots + a_n, a_i \in A_i, i \in \{1, \ldots, n\}.$$

We see that $H' = \pi(H)$ is a hyperplane of $A'$ and $\mathrm{Ker}(\pi) = A_1$. As before, $H'$ is not a two-sided ideal of $A'$, and $A'$ does not satisfy (ii). Since $\dim_F A' < \dim_F A$, by inductive hypothesis, there exists $u' \in H' \cap U(A') \subseteq H$. Choose any $u_1 \in U(A_1)$. Then $u = u' + u_1 \in U(A) \cap H = \emptyset$, a contradiction.

Because $A_n \not\subseteq H$, we see that $H_n = H \cap A_n$ is a hyperplane in $A_n$. Because $A'' = A_1 \dotplus \ldots \dotplus A_{n-1} \not\subseteq H$, we see that $H'' = H \cap A''$ is a hyperplane in $A''$. If there exists $u_n \in H_n \cap U(A_n)$ and $u'' \in H'' \cap U(A'')$, then $u = u_n + u'' \in H \cap U(A) = \emptyset$, a contradiction. So either $H_n \cap U(A_n) = \emptyset$ or $H'' \cap U(A'') = \emptyset$.

Suppose that $H_n \cap U(A_n) = \emptyset$. Then $A_n \cong F$ by Lemma 7. Therefore $A = A'' \dotplus A_n \cong A'' \oplus F$. Since $A$ does not satisfy (ii), it follows that $A''$ does not satisfy (iii). Theorem 3 holds for $A''$ and $H''$ because $\dim_F A'' < \dim_F A$. Therefore, there exists $u'' \in U(A'') \setminus H''$. Since $\dim_F A/H = 1$ and $u'' \notin H$, we have that $A/H = (u'' + H)F$. Let $u_n \in A_n \setminus H_n$. Then $u_n \notin H$, so that $A/H = (u_n + H)F$. Since $u'' + H$ and $u_n + H$ are $F$-linearly dependent, we have that $u'' + \alpha u_n \in H$ for some $0 \neq \alpha \in F$. However, $A_n \cong F$ implies that $\alpha u_n \in U(A_n)$. Then $u'' + \alpha u_n \in U(A'') + U(A_n) = U(A)$. It follows that $U(A) \cap H \neq \emptyset$, a contradiction. We conclude that $H'' \cap U(A'') = \emptyset$.

Lemma 8 gives us some $u_n \in U(A_n) \setminus H_n$. Let $u'' \in U(A'')$. Because $u'', u_n \in A \setminus H$, there is $0 \neq \alpha \in F$ such that $u = u'' + \alpha u_n \in H$. Then $u \in U(A) \cap H = \emptyset$, a contradiction. We conclude that Theorem 2 holds for $A$ and $H$.

It remains to prove Theorem 3 for $A$ and $H$. Assume it is false. Then $U(A) \subseteq H$ and $A$ does not satisfy (iii). By Lemma 6, $\mathrm{nil}(A) \subseteq H$. As in the above proof, this implies that $J(A) = \{0\}$, so that $A$ is the internal direct sum $A_1 \dotplus \ldots \dotplus A_n$ of $n \in \mathbb{N}$ two-sided ideals $A_i$ each of which is a simple $F$-subalgebra of $A$. As in the above proof, $A_i \not\subseteq H$ for all $i \in \{1, \ldots, n\}$. Moreover, $n \geq 2$ by Lemma 8.

The intersections $H_n = H \cap A_n$ and $H'' = H \cap A''$ are now hyperplanes in $A_n$ and $A'' = A_1 \dotplus \ldots \dotplus A_{n-1}$, respectively. Lemma 8 gives us some $u_n \in U(A_n) \setminus H_n$, and hence $u_n \notin H$. If there is some $u'' \in H'' \cap U(A'')$, then $u = u_n + u'' \in U(A) \setminus H$. This contradicts our assumption that $U(A) \subseteq H$. Hence $H'' \cap U(A'') = \emptyset$. Because $A$ does not satisfy (iii), its epimorphic image $A''$ does not satisfy (ii). So Theorem 2, for the algebra $A''$ with $\dim_F A'' < \dim_F A$, tells us that $H''$ is a two-sided ideal of $A''$. Since $\dim_F A''/H'' = 1$, we conclude that one of the direct summands $A_1, \ldots, A_{n-1}$ of $A''$, say $A_{n-1}$, is isomorphic to $F$, and that $H'' = A_1 \dotplus \ldots \dotplus A_{n-2}$. But no $A_i$ is contained in $H$ and $H'' \subseteq H$. Therefore $H'' = \{0\}$, $n = 2$, and $A'' = A_{n-1} = A_1 \cong F$.

We claim that $A_n \cong F$ and $H_n = \{0\}$. If there exists $u_n \in H_n \cap U(A_n)$, then $u'' + u_n \in U(A) \setminus H$ for any $0 \neq u'' \in A''$. This contradicts our assumption that $U(A) \subseteq H$. Hence $H_n \cap U(A_n) = \emptyset$. Then Lemma 7 implies that $A_n \cong F$ and $H_n = \{0\}$.

Now $A = A'' \dotplus A_n \cong F \oplus F$. If $|F| \geq 3$, then no one-dimensional subspace of $F \oplus F$ can contain $U(F \oplus F)$. Hence $|F| = 2$, so that $A \cong \mathbb{F}_2 \oplus \mathbb{F}_2$. Thus $A$ satisfies (iii), a contradiction. This contradiction shows that $A$ must satisfy Theorems 2 and 3. So the simultaneous proof of these theorems is complete.

## 4 Main Result

**PROOF OF THEOREM 1.** We note that once we prove that (i) $\iff$ (ii) $\iff$ (iii) $\iff$ (iv), then we can conclude that the groups $C/K^*$ and $\text{Hom}_{\mathbb{Z}}(G, K^*)$ are isomorphic via the map $(\gamma K^*)(\tau) = \frac{\tau(\gamma)}{\gamma}, \gamma \in C$. This map is surjective because the group $\mathrm{H}^1(G, L^*)$ is trivial (see, for example, [8]).

We see that (ii) $\implies$ (i).

Next we prove that (i) $\implies$ (ii). We begin by making the following observations.

(1) If $\gamma_1, \gamma_2 \in L$ are such that $\gamma_1 B, \gamma_2 B \subseteq B$, then $B \supseteq \gamma_1 B \supseteq \gamma_1 \gamma_2 B$.
(2) If $\alpha \in K^*$, then $\alpha B = B$.
(3) If $\tau \in G$, then $\tau(B) = B$.

Let $\gamma_1 = \gamma$ and $\gamma_2 = \prod_{1 \neq \tau \in G} \tau(\gamma)$. Then $\gamma_1 B, \gamma_2 B \subseteq B$, and at the same time, $\gamma_1 \gamma_2 B = B$ since $\gamma_1 \gamma_2 \in K^*$. Therefore, $B \supseteq \gamma_1 B \supseteq \gamma_1 \gamma_2 B = B$.

We have proved that (i) $\iff$ (ii).

We now show that (iii) $\implies$ (iv). Let $\gamma_\tau = \frac{\tau(\gamma)}{\gamma}$, then $\gamma_\tau \in K^*$. Let $S_\gamma = \{\gamma_\tau \mid \tau \in G\}$. Then $\gamma_{\tau_1} \gamma_{\tau_2} = \gamma_{\tau_1 \circ \tau_2}$ for all $\tau_1, \tau_2 \in G$. Therefore, $S_\gamma$ is a group under multiplication. It is a subgroup of the group

$$S = \{\mu \in K^* \mid \mu^{|G|} = 1\}$$

of order $w$. Hence, by Lagrange's Theorem, it would be sufficient to show that $\gamma^{|S_\gamma|} \in K^*$. This is indeed the case: $1 = \gamma_\tau^{|S_\gamma|} = \frac{\tau(\gamma^{|S_\gamma|})}{\gamma^{|S_\gamma|}}$ for any $\tau \in G$.

Next we show that (iv) $\implies$ (iii). From (iv) we conclude that $(\frac{\tau(\gamma)}{\gamma})^w = 1$ for any $\tau \in G$. We claim that $\frac{\tau(\gamma)}{\gamma} \in S$ (since $S \subseteq K^*$, we would be done). Indeed, each element of $S$ is a solution of the equation $x^w - 1 = 0$. This equation has at most $w$ solutions in any field extension of $K$. On the other hand, $|S| = w$, so that $S$ is the set of all solutions of this equation in any field extension of $K$. Because $(\frac{\tau(\gamma)}{\gamma})^w = 1$, we conclude that $\frac{\tau(\gamma)}{\gamma} \in S$.

We have proved that (iii) $\iff$ (iv).

We note that (iii) $\implies$ (i) because any element of $K$ is fixed by any element of $G$ and the map

$$\widehat{\gamma} : L \to L, t \mapsto t\gamma$$

is a bijection.

Finally, we show that (ii) $\implies$ (iii). This would prove that

$$(\text{i}) \iff (\text{ii}) \iff (\text{iii}) \iff (\text{iv}).$$

We know that $B \neq \emptyset$ by the Normal Basis Theorem (see, for example, [8]). Let $a \in B$. Define the map

$$\widetilde{a} : KG \to L, \sum_{\tau \in G} \alpha_\tau \tau \mapsto \sum_{\tau \in G} \alpha_\tau \tau(a).$$

Then

(1) $\widetilde{a}$ is an isomorphism of $KG$-modules.
(2) $\widetilde{a}(\omega KG) = N$.
(3) There exists a unique $K$-linear map $\Gamma$ which makes the following diagram commutative:

$$\begin{array}{ccc} L & \xrightarrow{\widehat{\gamma}} & L \\ \widetilde{a} \uparrow & & \widetilde{a} \uparrow \\ KG & \xrightarrow{\Gamma} & KG. \end{array}$$

We claim that $\gamma N$ is a $KG$-submodule of $L$. Because $\widetilde{a}^{-1}(\gamma N) = \Gamma(\omega KG)$, Theorem 4 would imply our claim if $\Gamma(\omega KG) \cap U(KG) = \emptyset$. Suppose that $\Gamma(\omega KG) \cap U(KG) \neq \emptyset$. Let $r \in \Gamma(\omega KG) \cap U(KG)$. Then

(1) $\{\tau r \mid \tau \in G\}$ is a $K$-basis of $KG$.
(2) $\widetilde{a}(r) = \gamma \widetilde{a}(r')$, where $r' \in \omega KG$ is such that $r = \Gamma(r')$.

Therefore, $\{\tau(\widetilde{a}(r)) \mid \tau \in G\}$ is a $K$-basis of $L$. On the other hand, since $\widetilde{a}(r') \in N$, a proper $KG$-submodule of $L$, we see that $\{\tau(\widetilde{a}(r')) \mid \tau \in G\}$ is not a $K$-basis of $L$. We conclude that $\gamma^{-1}B \not\subseteq B$, contradicting (ii). Hence, we have proved that $\gamma N$ is a $KG$-submodule of $L$.

Because $N, \gamma N$ are $KG$-submodules of $L$, we conclude that, for any $\tau \in G$,

$$\gamma N = \tau(\gamma N) = \tau(\gamma)\tau(N) = \tau(\gamma)N.$$

Therefore, $N$ is a vector space over the field $K' = K(\frac{\tau(\gamma)}{\gamma})$. Consequently, $\dim_K K'$ divides $\dim_K N = \dim_K L - 1$. Because $L$ is a $K'$-vector space, we conclude that $\dim_K K'$ divides $\dim_K L$. Hence $\dim_K K'$ divides $\gcd(\dim_K L - 1, \dim_K L) = 1$, i.e. $K(\frac{\tau(\gamma)}{\gamma}) = K$ for any $\tau \in G$. The proof is complete.